\documentclass[11pt]{article}
\usepackage{latexsym,amssymb,amsmath}
\def\cA{{\mathcal A}}
\def\cB{{\mathcal B}}
\def\cC{{\mathcal C}}
\def\cP{{\mathcal P}}
\def\cQ{{\mathcal Q}}
\def\cM{{\mathcal M}}
\def\cI{{\mathcal I}}
\def\cJ{{\mathcal J}}

\def\cF{{\mathcal F}}

\def\cY{{\mathcal Y}}
\def\cZ{{\mathcal Z}}
\def\cO{{\mathcal O}}
\def\cH{{\mathcal H}}

\def\CC{\mathbb C}
\def\RR{\mathbb R}
\def\HH{\mathbb H}
\def\AA{{\mathbb A}}

\def\OO{\mathbb O}

\def\ZZ{\mathbb Z}

\def\11{\mathbf 1}
\def\PP{\mathbb P}

\def\e1{\varepsilon_1}
\def\e2{\varepsilon_2}
\def\e3{\varepsilon_3}

\def\P2{{\PP}^2}

\def\00{\underline{0}}
\def\J0{{\cal J}_3(\underline{0})}

\def\PJ0{\PP({\cal J}_3(\underline{0}))}

\def\l{\lambda}
\def\a{\alpha}

\def\s{\sigma}

\def\e{\varepsilon}
\def\ot{{\mathord{\,\otimes }\,}}
\def\op{{\mathord{\,\oplus }\,}}

\def\lra{{\mathord{\;\longrightarrow\;}}}
\def\ra{{\mathord{\;\rightarrow\;}}}

\def\we{{\mathord{{\scriptstyle \wedge}}}}

\def\AP2{{\AA\PP}^2}
\def\RP2{{\RR\PP}^2}
\def\CP2{{\CC\PP}^2}
\def\HP2{{\HH\PP}^2}
\def\OP2{{\OO\PP}^2}

\newcommand\proof{{\noindent {\em Proof}.}\hspace{2mm}}
\newcommand\qed{{\hfill\hfill $\Box$}}

\newtheorem{theo}{Theorem}
\newtheorem{coro}[theo]{Corollary}
\newtheorem{lemm}[theo]{Lemma}
\newtheorem{prop}[theo]{Proposition}

\begin{document}

\title{Cubic hypersurfaces and integrable systems}
\author{Atanas Iliev, \ Laurent Manivel}


\maketitle
\footnotetext[1]{\ MSC 14J40, 14J45, 14K30, 14D06}
\footnotetext[2]{\ The 1-st author is partially supported by grant MI1503/2005
                  of the Bulgarian Foundation for Scientific Research}
                  
\begin{abstract} Together with the cubic and quartic threefolds, the cubic fivefolds 
are the only hypersurfaces of odd dimension $> 1$ for which the intermediate
jacobian is a nonzero principally polarized abelian variety
 (p.p.a.v.)). 
In this paper we show that 
the family of $21$-dimensional intermediate jacobians of cubic fivefolds containing
a given cubic fourfold $X$ is generically an algebraic integrable system. 
In the proof we apply an integrability criterion, introduced and used by
Donagi and Markman to find a similar integrable system over the family of cubic 
threefolds in $X$. To enter in the conditions of this criterion, we write down
explicitly the known by Beauville and Donagi symplectic structure on the family
$F(X)$ of lines on the general cubic fourfold $X$, and prove that the family of 
planes on a cubic fivefold containing $X$ is embedded as a Lagrangian surface in 
$F(X)$. By a symplectic reduction we deduce that our integrable system induces 
on the nodal boundary another integrable system, interpreted generically as the 
family of $20$-dimensional intermediate jacobians of Fano threefolds of genus 
four contained in $X$. Along the way we prove an Abel-Jacobi type isomorphism 
for the Fano surface of conics in the general Fano threefold of genus 4, and 
compute the numerical invariants of this surface.  
%
%
\end{abstract}

\section{Introduction}

\subsection*{1.1 Background} 

For a smooth compact K\"ahler manifold 
$Y$ of dimension $n$ and any positive integer $q \le n$, its 
$q$-th intermediate jacobian 
$$
J_q(Y) = (H^{q-1,q}(Y) \oplus\cdots\oplus H^{0,2q-1}(Y))/H^{2q-1}(Y,\ZZ)
$$
is a complex torus. Because of the skew-symmetry of the Riemann-Hodge 
bilinear relations the tori $J_q(Y)$ aren't in general 
abelian varieties, except for the Picard variety
$J_1(Y) = Pic^0(Y)$  
and the Albanese variety $J_n(Y) = Alb(Y)$, see Ch.2 \S 6 in \cite{CG}. 
If $n = dim\ Y = 2q-1$ is odd then $J(Y)=J_q(Y)$ is simply called  
the {\it intermediate jacobian} of $Y$. 
In the particular case when $Y \subset \PP^{n+1}$ 
is a smooth hypersurface of degree $d$, 
the intermediate jacobian $J(Y)$ happens to be a non trivial abelian 
variety only when $n = 1$ and $d \ge 3$,  
$n = 3$ and $d = 3, 4$ 
(i.e. when $Y$ is a cubic or quartic 3-fold),
or $n = 5$ and $d = 3$  (i.e. when $Y$ is a cubic $5$-fold), 
see p.43 in \cite{Col}.  

%
%

In \cite{DM1} Donagi and Markman studied certain complex algebraic 
analogs of integrable systems in mechanics. These analogs take as models the 
integrable Hamiltonian systems $f: M \rightarrow B$, 
where the phase space $M$ is a real symplectic $2n$-fold, 
and the general fibers of $f$ are compact Lagrangian $n$-folds in $M$, 
which by the Arnold-Liouville theorem are $n$-fold real tori.
In \S 7 of \cite{DM1} are given general criteria 
(the cubic conditions) which ensure when an algebraic family 
$f: M \rightarrow B$ of complex tori
represents an algebraic integrable system;  
for more details see \S 2 in \cite{DM1} or \S 2 in \cite{Fr}.   

The main application of these criteria  
is the observation that the relative intermediate jacobians 
over the gauged moduli spaces of Calabi-Yau threefolds are
in fact such algebraic integrable systems, see \cite{DM2}.  
This gives rise to the question whether besides the Calabi-Yau 
integrable systems there exist other nontrivial examples of 
families of intermediate jacobians that are also algebraic 
integrable systems, see e.g. Question 4.6 in \cite{AS}. 
It seems that until now, the only known such 
nontrivial example (besides the integrable systems coming 
from families of jacobians of curves) is given in Example 8.22 
of \cite{DM1}:

\smallskip

{\it If $X \subset \PP^5$ is a smooth cubic fourfold, and 
$B = \cH_s$ is the open subset in the dual projective space 
 parameterizing the smooth hyperplane sections 
$Y_b$ of $X$, then the relative intermediate jacobian 
$f: \cJ(\cY) \rightarrow \cH_s$ is an algebraic integrable system}. 
 
\smallskip 

This example is based on an abstract procedure, 
described in \S 8 of \cite{DM1}, 
that generates algebraic integrable systems in the context of
Lagrangian deformations. 
More precisely, for a given complex symplectic variety $W$, 
its Lagrangian Hilbert scheme parameterizes 
the Lagrangian subvarieties of $W$. By a result of Voisin and Ran 
the deformations of a smooth Lagrangian subvariety $F \subset W$ are 
unobstructed and the component $\cH$ of the Lagrangian Hilbert 
scheme of $W$ containing $F$ is generically smooth, see \cite{Voi, Ran}.
Consider the universal family $\cF \rightarrow \cH$ 
and its Picard bundle $Pic^o \cF \rightarrow \cH$. In this context,
Theorem 8.1 in \cite{DM1} states:

\smallskip

{\bf (DM)} \   
{\it   If $W$ is a smooth complex symplectic variety 
and $F$ is a smooth Lagrangian subvariety of $X$, 
then the symplectic structure on $W$ generates  
a natural symplectic structure 
on the relative Picard $Pic^o \cF \rightarrow \cH_s$  
over the base $\cH_s$ of smooth Lagrangian deformations of $F$ in $W$,
making the relative Picard $Pic^o \cF \rightarrow \cH_s$ 
an algebraic integrable Hamiltonian system.}

\smallskip

In the cited Example 8.22, the above criterion {\bf (DM)}  is applied 
in the case where $W = F(X)$ is the family of lines on a smooth cubic 
hypersurface $X \subset \PP^5$ and $F = F(Y) \subset F(X)$ is the Fano 
surface of lines on a smooth hyperplane section $Y$ of $X$. 
To enter in the conditions of {\bf (DM)}, the following results are used:

\smallskip
 
{\bf 1}. {\it The family $F(X)$ is a symplectic 4-fold}, see \cite {BD},    

{\bf 2}. {\it $F(Y) \subset F(X)$ is a Lagrangian surface in $F(X)$}, 
         see Ex.7 in \S 3 of \cite{Voi},  


{\bf 3}. {\it For a smooth cubic 3-fold $Y$, its intermediate jacobian 
$J(Y)$ is a principally polarized abelian variety of dimension $5$,  
the family $F(Y)$ of lines on $Y$ is a smooth surface of irregularity $5$, 
and the Abel-Jacobi map defines an isomorphism between 
$Alb F(Y)$ and $J(Y)$}, see \cite{CG}. 
  

\smallskip 
The most common situation when {\bf (DM)} can be applied 
is the case when the symplectic variety $W = S$ is a K3 surface 
and $F = C$ is a smooth curve of genus $g$ on $S$;
notice that any curve on a K3 surface $S$ is a Lagrangian subvariety
of $S$. In this lowest dimensional case, the base $\cH_s$ of the smooth 
Lagrangian deformation of $C$ in $S$ is an open subset of the complete 
linear system $|\cO_S(C)| \cong \PP^g$, and {\bf (DM)} yields 
that the relative jacobian $\cJ$ is a Lagrangian fibration over $\cH_s$.
When $S = S_{2g-2}$ is a K3 surface with a primitive 
polarization of genus $g = g(C)$ this fibration can be extended to 
a Lagrangian fibration $\overline{\cJ} \rightarrow \PP^g$ over the 
compactified 
relative jacobian $\overline{\cJ}$, which is a smooth compact 
complex symplectic variety birational to the $g$-th punctual Hilbert scheme
$Hilb_g S$, see \S 2 in \cite{B2}.
In this trend, J. Sawon, in a collaboration with K. Yoshioka,
has recently shown that for any $g,m \ge 2$ the $g$-th Hilbert power $Hilb_g T$
of the general primitive K3 surface $T$ of degree $m^2(2g-2)$ 
can be represented as a torsor over the compactified relative 
jacobian $\overline{\cJ} \rightarrow \PP^g$ of a primitive K3 
surface $S = S(T)$ of degree $2g-2$, thus proving the existence 
of a regular Lagrangian fibration over $\PP^g$ on the smooth compact 
complex sympectic variety $Hilb_g T$, see \cite{Saw}.  

\smallskip
The above examples give rise to the question when 
the relative Picard fibrations from {\bf (DM)} can be extended 
to Lagrangian fibrations including fibers over singular
Lagrangian deformations of $F \subset W$.  
%
For general integrable systems of type relative Picard as in {\bf (DM)},  
a partial compactification is described
in Theorem 8.18 of \cite{DM1} but the proof remains 
unpublished.

At the end of the same Example 8.22 of \cite {DM1} 
it is shown that the symplectic structure 
on $Pic^o({\cal F}) \rightarrow \cH_s$ can be extended to 
a still non-degenerate symplectic structure over the set 
$\cH_n \supset \cH_s$ of hyperplane sections $Y$ of $X$ that are 
allowed to be singular but with at most one node.
For a general nodal cubic 3-fold $Y_b \in \partial\cH_n:=\cH_n-\cH_s$
the family of lines $\ell \subset Y_b$ that pass through the node 
of $Y$ is parameterized by a smooth curve $C_b$ of genus $4$, 
and the generalized intermediate jacobian $J(Y)$ of $Y$ 
is a $\CC^*$-extension of the jacobian $J(C)$,
see e.g. \cite{CM}. Then the relative 
jacobian fibration $\cJ(\cY) \rightarrow \cH_n$
induces over the nodal boundary $\partial \cH_n$ a boundary 
fibration $\cJ(\cC) \rightarrow \partial \cH_n$ whose fibers 
are the jacobians $J(C_b)$ of the genus $4$ curves $C_b$. 
All this makes it possible to conclude that the boundary 
abelian fibration $\cJ(\cC) \rightarrow \partial \cH_n$ 
of the algebraic integrable system $\cJ(\cY) \rightarrow \cH_n$
is also an algebraic integrable system, at least over an open subset
of $\partial \cH_n$. Donagi and Markman call this system the 
{\it boundary integrable system} for $\cJ(\cY)$, 
see the end of \S 8 in \cite{DM1}. 
As communicated to us by Ron Donagi, this
is an instance of the algebraic symplectic reduction,
described in \S 2 of his later paper \cite{DP} with E. Previato.


\subsection*{1.2 Summary of the results in the paper}

In this paper we describe a new example of an algebraic 
integrable family of intermediate jacobians in the 
context of the integrability conditions {\bf (DM)}.    
It is an analog of the Example 8.22 from \cite{DM1}, where the 
fibers of the integrable system are the 5-dimensional intermediate 
jacobians of smooth cubic 3-folds contained in a fixed cubic fourfold as 
hyperplane sections. In our case the fibers of the integrable system 
are the 21-dimensional intermediate jacobians of the general cubic 
fivefolds containing the same cubic fourfold $X$ as a hyperplane section.    
By Theorem 13, the first main conclusion in our paper:

\smallskip
 
{\it The relative intermediate jacobian is an algebraic integrable system 
over an open subset of the family of cubic fivefolds containing a fixed 
general cubic fourfold $X$ as a hyperplane section}.  

\smallskip

Notice once again that the cubic fivefolds are the {\it unique}
hypersurfaces of odd dimension $> 3$ for which the intermediate 
jacobian is a non trivial abelian variety. 
Next, we study the degeneration of this integrable system on the 
boundary parameterizing the nodal cubic 
fivefolds through $X$, and prove that the induced boundary 
abelian fibration is also an algebraic integrable system.  

\smallskip

In our case the symplectic variety $W$ from {\bf (DM)}  is still the 
same as in the example of Donagi and Markman:
the 4-fold family $W = F(X)$ of lines on a fixed smooth 
cubic fourfold $X$. 
The difference is in the choice of the Lagrangian subvariety 
of departure $F \subset W$. 
At this point we should mention that discovering principally 
new or nontrivial examples of Lagrangian subvarieties of 
symplectic varieties, especially in a projective-geometric context,
looks like a happy occurrence. 
In our case this occurrence is realized as the Fano surface 
$F_2(Z)$ of planes on the general cubic fivefold $Z$ that
contains $X$ as a hyperplane section; the surface $F_2(Z)$
is regarded as a subvariety of $F(X)$ by the embedding given 
by the intersection-map $j_Z: F_2(Y) \rightarrow F(X)$, 
$\PP^2 \mapsto \PP^2 \cap X$.
We prove:

\smallskip

{\it The intersection image $j_Z(F_2(Z))$ of the
Fano surface $F_2(Z)$ of planes on the general cubic fivefold $Z$ containing
the cubic fourfold $X$ as a hyperplane section is a smooth Lagrangian 
surface inside the 4-fold family $F(X)$ of lines in $X$}, see Proposition 4. 

\smallskip






This relies on an explicit description of the symplectic form on 
$F(X)$, which is provided in \S 2.1.
Recall that in \cite{BD} the symplectic form on $F(X)$ 
is described only for the general Pfaffian cubic 4-folds $X$,
which form a divisor in the space of all cubic 
fourfolds $X \subset \PP^5$. 
For a general Pfaffian cubic fourfold $X$, the family $F(X)$ is known 
to be  isomorphic to the Hilbert square $Hilb_2 S$ of a K3 surface $S$ 
of genus $8$, and the symplectic form on $F(X)$ 
is the symplectic form on $Hilb_2 S$ described earlier by 
Fujiki and Beauville, see \cite{B1}.  
The existence of a symplectic form 
on $F(X)$ for the general cubic fourfold $X$ then follows 
by a deformation argument, see \cite{BD}. 

\smallskip

  
By Proposition 4 we enter in the conditions of {\bf (DM)}, 
this time with the symplectic fourfold $F(X)$ and its  
Lagrangian surface $j_Z(F_2(Z))$. This gives rise 
to a Lagrangian fibration on the relative Picard 
$Pic^o \cF \rightarrow \cH_s$ over the scheme $\cH_o$ 
of smooth deformations of $j_Z(F_2(Z))$ in $F(X)$, 
see the beginning of Section 5.  
In this situation, the analog of {\bf 3} above 
is a result of Collino: 

\smallskip

{\it The family $F_2(Z)$ of lines 
on the general cubic 5-fold $Z$ is a smooth surface of irregularity 
$21$, and the Abel-Jacobi map induces an isomorphim between the Albanese 
variety $Alb F_2(Z)$ and the intermediate jacobian $J(Z)$}, 
see \cite{Col}.

\smallskip

In fact $h^{5,0}(Z)$ and $h^{4,1}(Z)$ vanish, while 
$h^{2,1}(Z) = 21$, which explains why $J(Z)$ is a p.p.a.v. of 
dimension $21$.  By duality, $J(Z)$ is also isomorphic with 
the Picard variety $Pic^o F_2(Z)$.   This implies Theorem 13.


\medskip

Next, we study a partial compactification of the relative jacobian 
fibration from Theorem 13, by including to its base 
the codimension one boundary $\partial \cH_n$ 
of cubic fivefolds with one node.  
%
%
For a general nodal cubic 5-fold $Z_b$, the family 
of lines $\ell \subset Z_b$ which pass through the node of $Z_b$ 
is parameterized by a general smooth prime Fano threefold $Y_b$ 
of genus $4$. This threefold $Y_b$ is the analog of the genus 4 
curve $C_b$ from the boundary 
system in Example 8.22 of \cite{DM1}. 
 Since in this  situation 
the generalized intermediate jacobian $J(Z_b)$ is 
a $\CC^*$-extension of $J(Y_b)$, 
the boundary fibration $\cJ(\cY)$ has for  base the set 
of nodal cubic 5-folds $Z_b$ containing $X$ and 
as fibers the $20$-dimensional intermediate jacobians 
of their associated Fano 3-folds $Y_b$ of genus $4$. 
Such a Fano 3-fold $Y_b$ is a complete intersection 
of a quadric and a cubic in $\PP^5$. 
While the cubic can be chosen to be $X$, 
the quadric (identified with the base of the 
projective tangent cone to the node of $Z_b$) can move 
together with $b$, see \S 4.2. 

In \S 2 we find the analogs of {\bf 2} and {\bf 3} 
(i.e. boundary versions of Proposition 4 and Collino's 
Abel-Jacobi isomorphism) over the boundary $\partial \cH_n$ 
of nodal cubic fivefolds. 
To find these analogs we need for a given general nodal cubic 
fivefold $Z_b$ to express the Fano surface $F_2(Z_b)$ with its embedding
$j_Z$ in the family of lines $F(X)$ of the fixed cubic 4-fold $X$.
As shown in \S 2.2, 
the Fano surface $F_2(Z_b)$ of planes on $Z_b$ is almost 
the same as the Fano surface $F(Y_b)$ of conics on $Y_b$. 
More precisely, there is a natural map $j: F(Y_b) \rightarrow  F(X)$,
such that the image $j(F(Y_b))$ coincides with the isomorphic intersection
image of $F_2(Z_b)$ in $F(X)$. However while $F(Y_b)$ is smooth, 
the Fano surface $F_2(Z_b)$ has a double curve 
$\Gamma$, isomorphic to the family of lines on $Y_b$,
see Proposition 7. Nevertheless, we conclude 
that $j(F(Y_b)) = j_Z(F_2(Z_b))$ is a singular Lagrangian 
surface in $F(X)$, see Proposition 8.    

\smallskip

In \S 3.1 we study in more detail the Fano surface 
$F(Y)$ of the general prime Fano 3-fold $Y$ of 
genus $4$. In particular we find the invariants of $F(Y)$, 
see Corollary 10. 
%
Next, in \S 3.2 we prove the following analog of  Collino's 
Abel-Jacobi isomorphism for planes on cubic fivefolds:  

\smallskip

{\it The family $F(Y)$ of conics on the general prime Fano 3-fold $Y$ 
of genus $4$ is a smooth surface of irregularity $20$, 
and the Abel-Jacobi map induces an isomorphism between $Alb F(Y)$
and the intermediate jacobian $J(Y)$}, see Theorem 12. 

\smallskip
 
To prove Theorem 12, we follow the same program as Clemens and
Griffiths in their proof of the Abel-Jacobi isomorphism for the cubic
threefold. This program was subsequently applied in 
\cite{Let} to prove a similar Abel-Jacobi isomorphism for the Fano 
surface of conics on the general quartic hypersurface in $\PP^4$, then
in \cite{Col} 
for the Fano surface 
of planes on the general cubic 5-fold. 
In brief, this program 
consists in verifying that in a general Leftschetz pencil 
$\{ Y_t : t \in \PP^1 \}$ of Fano 3-folds of genus $4$, 
for a finite number of values of $t$ 
either the Fano surface $F(Y_t)$ acquires isolated singular points
but $Y_t$ remains smooth, or when $Y_t$ becomes singular then 
its singularity is a simple node and in this case $F(Y_t)$ 
has a smooth double curve for singular locus, 
see \cite{Let} and \S 3.2. 
As above, the Abel-Jacobi isomorphism
makes it possible to identify $Pic^o F(Y)$ with $J(Y)$.

In \S 4.2 we use these results together with the algebraic 
symplectic reduction procedure from \cite{DP}
to get Theorem 17, the second main conclusion 
of this paper:

\smallskip

{\it  The algebraic integrable system from Theorem 13 induces 
on the nodal boundary $\partial \cH_n$ a fibration by intermediate 
jacobians of Fano 3-folds of genus $4$ which is generically 
an algebraic integrable system.}

\medskip

\noindent {\it Acknowledgements}.
We express our gratitude to Ron Donagi for pointing us to 
his paper with E. Previato, 
to Vesselin Drensky for the assistance in parallel computing the 
invariants from \S 3 with Maple.
The first author uses the occasion to thank after many years 
Alessandro Verra, Giuseppe Ceresa and especially Maurizio Letizia 
for the interesting conversations around the Clemens-Griffiths program.  


\section{Revisiting the family of lines on the cubic fourfold}

\subsection*{2.1 The symplectic form}


Let $X\subset\PP V=\PP^5$ be a general cubic hypersurface, 
with equation $P=0$ for some general polynomial $P\in S^3V^*$. 
The Fano variety  $F(X)$ of lines contained in $X$ is a subvariety 
of the Grassmannian $G(2,V)$. It can be defined as the zero locus 
of the global section $s_P$ of the vector bundle $S^3T^*$, naturally
defined by $P$, where $T$ denotes the rank two tautological bundle
on the Grassmannian.  Since $s_P$ is a general section of the 
globally generated vector bundle $S^3T^*$, $F(X)$ is a smooth 
four dimensional variety. 

It was proved by Beauville and Donagi that $F(X)$ has a 
symplectic structure. Indeed they showed that for $X$ a Pfaffian cubic 
hypersurface there exists a K3 surface $S$ such that $F(X)=Hilb_2S$,
which is well known to inherit a symplectic structure from that of $S$. 
Then the general case follows by a deformation argument, see \cite{BD}.
Nevertheless
the use of deformations makes that the symplectic form on $F(X)$ is 
not explicit. Markushevich and Tikhomirov \cite{MT} showed how to 
deduce it from the Voisin's observation that the Fano varieties of 
lines in the hyperplane sections of $X$ are Lagrangian surfaces 
in $F(X)$ (see section 3). We shall give below a completely explicit 
and down-to-earth description of this form. 
In fact it suffices to exhibit a non zero holomorphic two form on 
$F(X)$ -- indeed since $h^{2,0}(F(X))=1$ 
we know that there exists a unique such form, 
so it must define the symplectic structure, that it, it will
automatically be a non degenerate closed form.  
\smallskip
Consider the  tangent exact sequence 
$$0\ra TF(X)\ra TG(2,V)_{|F(X)}\ra S^3T^*_{|F(X)}\ra 0.$$
Remember that $TG(2,V)\simeq Hom(T,Q)$, where $Q=V/T$ denote the 
quotient bundle on the Grassmannian. Therefore, the tangent space
to $TF(X)$ at a line $\ell=\PP T$ is
$$T_{\ell}F(X)=\{u\in Hom(T,Q), \quad P(x,x,\bar{u}(x))=0\;
\forall x\in T\},$$
where $P(x,y,z)$ is the polarization of $P(x)$ and $\bar{u}(x)\in V$
is any representative of $u(x)\in Q$. 

Note that $P(x,x,.)=0$ is an equation of $T_xX$. If $\ell$ is a general
line on $X$ and $x$ varies in $\ell$, we get a quadratic pencil of
hyperplanes in $\PP V$, whose intersection is a plane $\Pi=\PP S
\supset\ell=\PP T$. In particular, $T_{\ell}F(X)\supset Hom(T,S/T).$  

There is a natural skew-symmetric form on $T_{\ell}F(X)$, defined 
up to scalar as follows. Choose some vectors $e,f,g,h$ of $T$ and let, for 
$u,v\in T_{\ell}F(X)$,
$$\Omega(u,v;e,f,g,h)=
P(e,e,u(g))P(f,f,v(h))-P(e,e,v(g))P(f,f,u(h)).$$
This can be seen as defining a skew-symmetric bilinear map on
$S^2T\otimes T$ with values in $\wedge^2T_{\ell}F(X)^*$. 
Note that $S^2T\otimes T=S^3T\oplus T\ot\wedge^2T      $. We get the component 
$S^3T$ by letting $g=e$ in the formula above, which gives zero on 
$T_{\ell}F(X)^*$. So we need only consider the other component,
which gives a map $(\wedge^2T)^3\ra \wedge^2T_{\ell}F(X)^*$ 
defined by the formula
$$\begin{array}{rcl}
\omega(u,v;e,f) &= & P(e,e,u(f))P(f,f,v(e))-P(e,e,v(f))P(f,f,u(e)) \\
 & & +2P(e,f,u(f))P(e,e,v(f))-2P(e,e,u(f))P(e,f,v(f)) \\
 & & +2P(f,f,u(e))P(e,f,v(e))-2P(e,f,u(e))P(f,f,v(e)).
\end{array}$$
We can see $\omega$ as a form with values in the line bundle
$(\wedge^2T^*)^{\otimes 3}=\cO(3)$. 
This form has rank two at the generic point of $F(X)$ and its radical
is precisely $Hom(T,S/T)$. 

There is also a natural quadratic form on $\wedge^2T_{\ell}G(2,V)$, 
defined again up to scalar by the following formula: if $u,v,u',v'\in
Hom(T,Q)$ and $e,f$ is a basis of $T$, let 
\begin{eqnarray*}
K(u\we v,u'\we v')= & u(e)\we u'(f)\we v(e)\we v'(f)
-u(f)\we u'(f)\we v(e)\we v'(e) \\
  & +u(f)\we u'(e)\we v(f)\we v'(e)
-u(e)\we u'(e)\we v(f)\we v'(f),
\end{eqnarray*}
seen as a element of the line $\we^4Q$. 
More precisely, $K$ is a quadratic form on $\wedge^2T_{\ell}G(2,V)$
with values in the line bundle $Hom((\wedge^2T),\wedge^4Q)=\cO(3)\ot\det(V)$. 

Observe that this form can be described in terms of the natural
decomposition $\wedge^2T_{\ell}G(2,V)=\wedge^2T^*\ot S^2Q\op
S^2T^*\ot \wedge^2Q$, as the composition of the natural morphisms
$$S^2(\wedge^2T_{\ell}G(2,V))\ra S^2(S^2T^*\ot \wedge^2Q)
\ra S^2(S^2T^*)\ot S^2(\wedge^2Q)\ra (\wedge^2T^*)^{\ot 2}\ot\wedge^4Q.$$

Now we restrict the quadratic form $K$ to $\wedge^2T_{\ell}F(X)$. 
We claim that the restriction has rank at most two. Indeed, 
it is clear that $K$ vanishes on $Hom(T,S/T)\wedge T_{\ell}F(X)$,
a hyperplane in $\wedge^2T_{\ell}F(X)$. In particular we can
write 
$$
K(u\we v,u\we v) = \omega(u,v)\Omega(u,v)$$
for $u,v\in T_{\ell}F(X)$, where $\Omega$ is now a 
well-defined skew-symmetric form on $T_{\ell}F(X)$ (with values 
in $\det(V)$, to be precise), since
$K$ and $\omega$ both take values in $\cO(3)$. 

\begin{theo}
The form $\Omega$ defines a symplectic structure on the
Fano variety $F(X)$. 
\end{theo}

\proof  As explained above we just need to prove that $\Omega$
is non trivial at some point $\ell=\PP T$
of $F(X)$. This reduces to an easy computation in coordinates:
we chose a basis $e_1,\ldots ,e_6$ of $V$ such that $T=\langle e_1,
e_2\rangle$. If $z_1,\ldots ,z_6$ are the corresponding coordinates 
on $V$, we may suppose that $P(e_1,e_1,z)=z_4$, $P(e_1,e_2,z)=z_5$, 
$P(e_2,e_2,z)=z_6$. In particular the intersection of the tangent 
spaces of $X$ along $\ell$ is $S=\langle e_1,e_2,e_3\rangle$.

An element of the tangent space $T_{\ell}F(X)$ is a homomorphism
$u\in Hom(T,Q)$ such that 
$$\begin{array}{rcl}
u(e_1) & = & ae_3+be_5-2ce_6,\\
u(e_2) & = & de_3-2be_4+ce_5
\end{array}$$ 
for some scalars $a,b,c,d$. Take another $u'\in T_{\ell}F(X)$
defined by the scalars $a',b',c',d'$. Then $K(u\we u',u\we u')$
is proportional to $u(e_1)\we u(e_2)\we u'(e_1)\we u'(e_2)$, that 
is, to the determinant
$$\left\|\begin{array}{cccc}
 a & d & a' & d' \\ 0 & -2b & 0 & -2b' \\
 b & c & b' & c' \\  -2c & 0 & -2c' & 0 
\end{array}\right\|
=4(bc'-b'c)(a'c-ac'+bd'-b'd).$$
The factor $bc'-b'c$ is proportional to $\omega(u,u')$. 
We conclude that $\Omega(u,u')$ is proportional to $a'c-ac'+bd'-b'd$,
so that $\Omega$ has maximal rank.\qed 

\medskip Note that $\omega(u,u')=0$ if and only if the 
tangent plane 
generated by $u$ and $u'$ contains an element $u''\in Hom(T,Q)$
whose image is contained in $S/T$. If it is not the case, then 
this tangent plane is isotropic if and only if 
$$u(T)+u'(T)\neq Q,$$
that is, there is a hyperplane $H\subset Q$ containing the image 
of any $u''$ in the tangent plane.

\subsection*{2.2 Lagrangian surfaces in $F(X)$}

Our very explicit description of the symplectic structure of $F(X)$
will help us to identify Lagrangian surfaces. 

\subsubsection*{2.2.1 Voisin's example}
First we recover the following observation of  
Claire Voisin (see \cite{Voi}). 

\begin{prop} Let $Y=X\cap H$ be a generic hyperplane section. 
Then the Fano surface $F(Y)$ of lines in $Y$ is a Lagrangian 
surface in $F(X)$. 
\end{prop}

\proof For $u,u'$ in $T_{\ell}F(Y)
\subset Hom(T,Q)$, the images $u(T)$ and $u'(T)$ are contained 
in the hyperplane of $Q$ defined by $H$. So clearly $K(u\we u',u\we
u')=0$. But $\omega (u,u')$ is non zero, at least generically, 
so $\Omega(u,u')=0$.  \qed 

\subsubsection*{2.2.2 Planes in cubic fivefolds} 
We now turn to the inverse situation where $X$ is a hyperplane
section of a cubic fivefold $Z\subset\PP^6=\PP W$, say $X=Z\cap H$
with $H=\PP V$. Generically, such a cubic fivefold contains projective 
planes, and the Fano variety $F_2(Z)$ of planes in $Z$ is a smooth surface 
in the Grassmannian $G(3,W)$. Cutting such a plane $\Pi$ with the
hyperplane $H$ we get a line in $X$, because $X$ contains no plane.
Hence a map 
$$i_Z:\; F_2(Z)\lra F(X).$$

\begin{lemm}
For general $X$ and $Z$, the map $i_Z$ is a closed embedding. 
\end{lemm}

\proof Two planes in a general cubic fivefold cannot meet 
along a line, so $i_Z$ is injective. We check it is also immersive.

Let $\Pi=\PP S$ be a plane in $Z$, such that $\ell=P\cap H$ is a line
in $X$. If $R=0$ is an equation of $Z$, we have
$$T_{\Pi}F_2(Z)=\{u\in Hom(S,W/S),\quad R(x,x,u(x))=0\;
\forall x\in S\}.$$
The differential of $i_Z$ at $\Pi$ maps $u\in Hom(S,W/S)$ to its 
restriction $u_{|T}$ from $T$ to $V/T\simeq W/S$. If this restriction 
is zero, then $u=e^*\ot f$ has rank one and $R(x,x,f)=0$ for all 
$x\in S$. This means that the intersection of the tangent hyperplanes 
to $Z$ along $\Pi$ intersect along a linear space strictly larger 
that $\Pi$. Call such a plane {\it special}. We  claim that a general 
cubic fivefold contains no special plane. 
Indeed, choose linear coordinates $x_0,\ldots ,x_6$ on $\PP^6$, 
such that the plane $\Pi$ is defined by $x_3=x_4=x_5=x_6=0$. If 
$Z$ contains $\Pi$, its equation is of the form
$x_3Q_3+x_4Q_4+x_5Q_5+x_6Q_6$ for some quadrics $Q_3,\ldots ,Q_6$. 
The intersection of the tangent hyperplanes to $Z$ along $\Pi$
contains $\Pi^+$ of equations $x_4=x_5=x_6=0$ if and only if the 
quadric $Q_3=0$ contains $\Pi$. So $\Pi$ is special if and only 
if it is special in the sense of \cite{Col}. Collino proves that 
if $\Pi$ is contained in the smooth locus of $Z$, then it defines a 
smooth point of $F_2(Z)$ if and only if it is non special. In particular
a smooth $Z$ such that $F_2(Z)$ is smooth contains no special plane, 
and $i_Z$ is immersive. \qed

\begin{prop}\label{lag2}
The map $i_Z$ embeds $F_2(Z)$ in $F(X)$ as a Lagrangian surface.
\end{prop}

\proof 
If $\omega$ vanishes on $(i_Z)_*
T_{\Pi}F_2(Z)$, then $T_{\Pi}F_2(Z)$ must contain a morphism 
$u\in Hom(S,W/S)$ whose restriction to $T$ has rank one. But 
this cannot happen generically, so we just need to prove that 
our quadratic form $K$ on 
$\wedge^2T_{\ell}F(X)$ vanishes on the line $\wedge^2(i_Z)_*
T_{\Pi}F_2(Z)$.

By definition of $K$ this means that two morphisms in 
$T_{\Pi}F_2(Z)$ cannot send $T$ to two independent subspaces 
in $W/S$. So our claim follows from the 
following observation:

\begin{lemm}
There exists a hyperplane $h_{\Pi}$ in $W/S$ such that 
for all $u\in T_{\Pi}F_2(Z)\subset Hom(S,W/S)$, we have
$u(S)\subset h_{\Pi}$. 
\end{lemm}

\proof Again we just need to prove this for a general cubic fivefold $Z$
and a general plane $\Pi\in F_2(Z)$. 
The cubic equation $R$ of $Z$ defines a linear 
map $\rho: W/S\ra Sym^2S^*$, mapping $\bar{w}$ to the quadratic
form $x\mapsto R(x,x,w)$ on $S$, for any representative 
$w$ of $\bar{w}$ in $W$. Consider the diagram
$$\begin{array}{ccc}
 & & \wedge^2S^*\ot S^* \\
 & & \downarrow \\
T_{\Pi}F_2(Z)\subset Hom(S,W/S) & \stackrel{id\ot\rho}{\lra} 
 & S^*\ot Sym^2S^* \\
 & & \downarrow \\
 & & Sym^3S^*
\end{array}$$
where the vertical strand is part of a Koszul complex. 
Since $T_{\Pi}F_2(Z)$ maps to zero in $Sym^3S^*$, its image 
by $id\ot\rho$ lies in the image of $\wedge^2S^*\ot S^*$. 
Note that since $S$ has dimension three, $\wedge^2S^*\ot S^*
=\det S\ot End(S)$. Once we chose a generator $\a$ of $\det S$,
we thus conclude that for any $u\in T_{\Pi}F_2(Z)$, there 
exists some $\theta_u\in End(S)$ such that 
$$R(x,x,u(y))\a = \theta_u(x)\we x\we y \qquad \forall x,y\in S.$$ 
Generically, the endomorphism $\theta_u$ is semisimple, and for any 
$z\in u(S)$ the quadratic form $x\mapsto R(x,x,z)$ vanishes 
along the three eigenlines of $\theta_u$. 

We interpret this as follows. The pull-back by $\rho$ of the 
discriminant defines a cubic surface $\Sigma$ in $\PP(W/S)$, and this 
surface is smooth by the genericity assumption. Since the 
image of the plane $\PP u(S)$ by $\rho$ is the net of conics 
passing through three general points, $\PP u(S)$ is a trisecant 
plane to the surface $\Sigma$. But there are only a finite number
of such trisecant planes, so $u(S)$ does not depend on $u$ (as long as
$u$ has maximal rank). This is what we wanted to prove. \qed

\subsubsection*{2.2.3 Conics on Fano threefolds of genus four}
A Fano threefold of genus four is the complete intersection $Y=Q\cap W$
of a quadric $Q$ and a cubic hypersurface $W$ in $\PP^5$. 

\begin{prop}
For $Y$ general, the set of lines in $Y$ is a smooth curve $\Gamma(Y)$
and the set of conics in $Y$ is a smooth surface $F(Y)$.
\end{prop}

\proof Choosing equations for $Q$ and $W$ we get sections 
$\s_Q$ and $\s_W$ of $S^2T^*$ and $S^3T^*$ on $G(2,6)$, and the 
zero locus of $\s_Q\op\s_W$, is precisely $\Gamma(Y)$. Hence
the first assertion. 

The Hilbert scheme $\cH$ of conics in $\PP^5$ is nothing else than 
the total space of the projective bundle $\PP(S^2T^*)\stackrel{\pi}{\ra}
G(3,6)$. On $\cH$ we have a tautological sequence 
$$0\ra \cO(-1)\ra\pi^*S^2T^*\ra\cQ\ra 0,$$ 
and the equation of $Q$ defines a section $\tau_Q$ of $\cQ$ whose zero-locus 
is  the set of conics contained in $Q$. 
In a similar way, the vector bundle $\cM$ defined by the exact sequence 
$$0\ra \pi^*T^*\ot\cO(-1)\ra\pi^*S^3T^*\ra\cM\ra 0,$$ 
has a global section $\tau_W$ defined by the equation of $W$
whose zero locus is the set of conics contained in $C$. 
Now $\tau_Q\op\tau_W$ is a general section of the globally generated 
vector bundle $\cQ\op\cM$, whose zero locus is $F(Y)$. So 
by Bertini, $F(Y)$ is a smooth surface. \qed

\medskip\noindent {\it Remark}.
The invariants of the curve $\Gamma(Y)$ were computed in 
\cite{Mar}, in particular its genus is $g(\Gamma(Y))=271$. We will
compute the invariants of $F(Y)$ in the next section.     
As noticed in Ch. 2 , \S 4 of \cite{Isk}, 
$F(Y)$ contains a special curve $B(Y) \subset F(Y)$ 
that consists of conics that lie in the ruling planes of $Q$. 
Since $Q$ has two family of ruling planes, this curve has two 
connected components $B_+(Y)$ and $B_-(Y)$. 

\medskip Now consider in $\PP^6$ a general nodal cubic 5-fold, 
$X$, with equation $x_0Q+W$, where $x_0$ is an equation of the 
hyperplane $\PP^5\supset Q,W$. Denote its node by $e_o = (1:0:...:0)$. 
The projective tangent cone $K(Q)$ to $X$ at $e_o$ 
is the cone in $P^6$ with vertex $e_o$ over the quadric 
$Q$.Its intersection with $X$ is the cone $K(Y)$ over the 
Fano 3-fold $Y = W \cap Q \subset P^5$. 

From another point of view, the cone $K(Y)$ is the union of 
all lines $\ell\subset X$ which pass through the node of $X$;
in other words, $Y$ is the base of the family of lines on $X$
that pass through $e_o$. 
Denote by $F_2(X)$ the Fano variety of planes in $X$.  
We shall construct a natural map 
$$
f: F(Y) \rightarrow F_2(X)
$$
as follows: 
Let $q \in Y$ be a conic. Then $q$, together with 
the node $e_o$ span a 3-space 
$P^3_q$, which intersects $X$ along a cubic surface 
$S_q = P^3_q \cap X$.
This surface is not irreducible, since it contains the 
quadratic cone $K(q)$ with vertex $e_o$ and base $q$. 
Therefore 
$$
S_q = K(q) + \PP^2_q,
$$ 
where $P^2_q\subset X$ is a plane, and we let 
$q \mapsto f(q) := P^2_q$.  

Note that $\ell_q:= P^2_q \cap W$ is a line 
that is 2-secant to $q$. We denote by $j$ the map $F(Y) \rightarrow
F(W)$ sending  $q$ to $\ell_q$. 

\medskip 
Next, we shall try to find where is defined the inverse of 
$f: F(Y) \rightarrow F_2(X)$.

First note that any plane $\Pi$ in $X$ that passes 
through $e_o$ evidently lies in $K(Y)$ and intersects on $Y$ 
a line $\ell=\Pi\cap Y$. Conversely, $\Pi$ is just the span 
of $\ell\cup e_0$. 
In other words, the curve of 
lines $\Gamma(Y)$ is the base of the family $\Gamma(X)$ of 
planes on $X$ that pass through $e_o$. 
in particular $\Gamma(X) \cong \Gamma(Y)$.

Let $\Pi\in F_2(X) - \Gamma(X)$, and let $\ell=\Pi\cap Y \in \Gamma(Y)$ be 
its corresponding line on $Y$. Since $e_o \notin\Pi$, the span 
of $\Pi\cup e_o$ is a 3-space, that intersects $X$ along a cubic surface 
$S_{\Pi}=\Pi\cup Q_{\Pi}$, where $Q_{\Pi}$ is a quadric surface in $X$. 
 
Since $S_{\Pi}$ clearly passes through $e_o = Sing(X)$, but $\Pi$
does not, the quadric $Q_{\Pi}$ must be singular at $e_o$.
Therefore $Q_{\Pi}$ must be a quadratic cone with vertex $e_o$, meeting 
$\PP^5$ along a conic $q$ on $Y = W \cap Q$. 
Clearly $q = f^{-1}(\Pi)$.

Thus $f$ is invertible outside $\Gamma(X) \cong \Gamma(Y)$. 
Their remains to find the pre-images on $F(Y)$ 
of the planes $\Pi\in \Gamma(X) \subset F_2(X)$.
Let $\ell=\Pi\cap Y \in \Gamma(Y)$. 

Let $q \in f^{-1}(\Pi)$. Then with the previous notations 
$S_q = P^3_q \cap X$ is a cone with vertex $e_0$, which can 
be defined as the set of points $(x_0,y)$, with $y$ in the span 
$\PP^2_q$ of $q$, such that $x_0Q(y)+W(y)=0$. So $Q$ must contain the
plane 
$\PP^2_q$, which itself contains $\ell$. But  a line in a 
four dimensional quadric is contained in exactly one plane in each 
of the two rulings of $Q$. Denote the two planes in $Q$ that contain
$\ell$ by $P_+$ and $P_-$. The intersection of $W$ with $P_{\pm}$
is the union of $\ell$ with a conic $q_{\pm}$, and we conclude that 
$$f^{-1}(\Pi)=\{q_+,q_-\}.$$
We summarize our discussion. 

\begin{prop}\label{sing}
The map $f: F(Y) \rightarrow F_2(X), \ q \mapsto P^2_q$
is an isomorphism outside the curve $B(Y)=B_+(Y)\cup B_-(Y)\
subset F(Y)$, and the restriction of $f$ to $B_{\pm}(Y)$ is 
an isomorphism with $\Gamma(Y)$. 

In particular $\Gamma(Y)$ is the singular locus of $F_2(X)$, 
whose normalization is $f$. 
\end{prop}

Now, since the condition of being Lagrangian is closed, we can 
deduce from Propositions \ref{lag2} and \ref{sing} the following result:

\begin{prop}
For any general cubic fourfold $W\subset\PP^5$ and any general
quadric $Q\subset\PP^5$, the map sending a conic $q\subset Y=Q\cap W$
to the residual line $\ell_q$ in the intersection of the linear 
span $P^2_q$ of $q$, with $W$, maps the surface $F(Y)$ of conics in $Y$
to a singular Lagrangian surface in $F(W)$. 
\end{prop}

\section{The family of conics on the Fano threefold of genus four}

\subsection*{3.1 Invariants of the Fano surface $F(Y)$}

Since $F(Y)\subset\cH$ is defined as the zero-locus of a general section 
of the vector bundle $\cQ\op\cM$, its fundamental class is given
by the Thom-Porteous formula, that is 
$$[F(Y)]=c_{12}(\cQ\op\cM)=c_5(\cQ)c_7(\cM)\in A^{12}(\cH).$$
Recall that the Chow ring $A(\cH)$ is generated over $A(G)$ by the 
hyperplane class $h=c_1(\cO(1))$, modulo the relation $c_6(\cQ)=0$. 
In particular $A(\cH)$ is a free $A(G)$-module with basis $1,\ldots
,h^5$. We say a class in $A(\cH)$ is written in normal form when it is 
expressed in that basis. We have
$$c(\cQ)=\frac{c(S^2T^*)}{1-h}, \qquad 
c(\cM)=\frac{c(S^3T^*)}{c(T^*(-1))}.
$$
In particular the first identity gives the relation
$hc_5(\cQ)=-c_6(S^2T^*)$.

The Chern classes of $S^2E$ and $S^3E$, for $E$ a vector bundle 
of rank three, are given by universal formulas in terms of the Chern
classes of $E$, that we denote by $c_1,c_2,c_3$. The splitting 
principle easily gives
$$c(S^2E)=(1+2c_1+4c_2+8c_3)(1+2c_1+c_1^2+c_2+c_1c_2-c_3).$$
The computation for $S^3E$ is much more complicated. With the help
of Maple we get the following formulas:
$$\begin{array}{ccl}
c_1(S^3E) &=& 10c_1 \\ 
c_2(S^3E) &=& 40c_1^2+15c_2 \\
c_3(S^3E) &=& 82c_1^3+111c_1c_2+27c_3 \\ 
c_4(S^3E) &=& 91c_1^4+315c_1^2c_2+189c_1c_3+63c_2^2\\ 
c_5(S^3E) &=& 52c_1^5+429c_1^3c_2+513c_1^2c_3+324c_1c_2^2+162c_2c_3\\ 
c_6(S^3E) &=& 12c_1^6+282c_1^4c_2+679c_1^3c_3+593c_1^2c_2^2+792c_1c_2c_3+
85c_2^3+27c_3^2\\ 
c_7(S^3E) &=& 72c_1^5c_2+448c_1^4c_3+464c_1^3c_2^2+1386c_1^2c_2c_3+
259c_1c_2^3+108c_1c_3^2+243c_2^2c_3\\ 
c_8(S^3E) &=& 120c_1^5c_3+132c_1^4c_2^2+1116c_1^3c_2c_3+
246c_1^2c_2^3+81c_1^2c_3^2+567c_1c_2^2c_3+36c_2^4+243c_2c_3^2\\ 
c_9(S^3E) &=& 360c_1^4c_2c_3+72c_1^3c_2^3+108c_1^3c_3^2+540c_1^2c_2^2c_3+
36c_1c_2^4-243c_1c_2c_3^2
+108c_2^3c_3+729c_3^3\\ 
c_{10}(S^3E) &=& 108c_1^4c_3^2+216c_1^3c_2^2c_3-486c_1^2c_2c_3^2+
108c_1c_2^3c_3+729c_1c_3^3 \\
\end{array}$$
If $x_1,x_2,x_3$ are the Chern roots of $E$, and if $L$ is a line bundle
with first Chern class $h$, we can also compute
$$\begin{array}{rcl}
c(E\ot L^*)^{-1} & = & \prod_i(1+x_i-h)^{-1}
=(1-h)^{-3}\prod_i(1+\frac{x_i}{1-h})^{-1} \\
 & = & \sum_{k\ge 0}(-1)^ks_k(E)(1-h)^{-k-3},
\end{array}$$
where $s_k$ denotes the $k$-th Segre class of $E$. 

Now we specialize to $E=\pi^*T^*$ (for simplicity we omit the 
symbol $\pi^*$ in the sequel) and $L=\cO(1)$, and we deduce that 
$$c_7(\cM)=\sum_{j,k\ge 0}(-1)^k
\binom{j+k+2}{j}c_{7-j-k}(S^3T^*)s_k(T^*)h^j.$$ 
We use the basis of $A(G)$ given by the Schubert classes $\s_{ijk}$, 
$3\ge i\ge j\ge k\ge 0$. Among these are the Chern classes 
$c_i(T^*)=\s_i$ and the Segre classes $s_k(T^*)=\s_{1^k}$
of $T^*$. In the Schubert basis the total Chern class of 
$S^2T^*$ is 
$$c(S^2T^*)=1+4\s_1+5\s_{2}+5\s_{11}+20\s_{3}+15\s_{21}
+30\s_{31}+10\s_{22}+6\s_{211}+20\s_{32}+12\s_{311}+4\s_{221}
+8\s_{321}.$$
Denote by $c_7(\cM)_{(j)}$ the coefficient of $h^j$ in that sum, 
for $0\le j\le 7$. Using the relation
$hc_5(\cQ)=-c_6(S^2T^*)=-8\s_{321}$, we get 
$$[F(Y)]=c_5(\cQ)c_7(\cM)_{(0)}-8\s_{321}\sum_{j=1}^7c_7(\cM)_{(j)}h^{j-1}.$$
This almost expresses $[F(Y)]$ in normal form, except for the 
term with $j=7$ in the sum. Since $c_7(\cM)_{(7)}=36$, 
the normal form is 
$$[F(Y)]=c_5(\cQ)c_7(\cM)_{(0)}-8\s_{321}\sum_{j=1}^6c_7(\cM)_{(j)}h^{j-1}
+288\s_{321}\sum_{i=0}^5c_{6-i}(S^2T^*)h^i.$$
The computations that remain are to be done  in $A(G)$, which is zero
in degree greater than nine. Since $[F(Y)]$ has degree $12$, it
must be a combination of 
the four classes $\s_{333}h^3$, $\s_{332}h^4$, $\s_{331}h^5$ 
and $\s_{322}h^5$.  

First we compute $c_7(\cM)_{(0)}=c_7(S^3T^*)-c_6(S^3T^*)\s_1+
c_5(S^3T^*)\s_{11}-c_4(S^3T^*)\s_{111}$. The formulas for the 
Chern classes of a third symmetric power give
$$\begin{array}{rcl}
c_7(S^3T^*) & = & 8820\s_{331}+6342\s_{322}, \\
c_6(S^3T^*) & = & 3675\s_{33}+7140\s_{321}+1302\s_{221}, \\
c_5(S^3T^*) & = & 2870\s_{32}+2436\s_{311}+1442\s_{211}, \\
c_4(S^3T^*) & = & 1155\s_{31}+560\s_{22}+588\s_{211}, 
\end{array}$$
hence $c_7(\cM)_{(0)}=1757\s_{331}+1190\s_{332}$. The other relevant 
classes for $\cM$ are 
$$\begin{array}{rcl}
c_7(\cM)_{(6)} & = & 196\s_1 \\
c_7(\cM)_{(5)} & = & 595\s_2+406\s_{11}\\
c_7(\cM)_{(4)} & = & 1375\s_3+1500\s_{21}+404\s_{111}
\end{array}$$
The final result of the computation is:

\begin{prop}\label{fund}
The fundamental class of the surface $F(Y)\subset\cH$ is 
$$[F(Y)]=15840\s_{333}h^3+8100\s_{332}h^4+(1341\s_{331}
+774\s_{322})h^5.$$
In particular, we have the following intersection numbers in the Chow
 ring of $F(Y)$:
$$h^2=0, \quad h\s_1=-360, \quad \s_2=1341, \quad \s_1^2=2115.$$
\end{prop}

Now we can compute the Chern numbers of $F(Y)$. From the exact sequences 
$$\begin{array}{c}
0\lra TF(Y)\lra T\cH_{|F}\lra \cQ_{|F}\oplus\cM_{|F}\lra 0, \\
0\lra T^v\cH\lra T\cH\lra \pi^*TG\lra 0 \\
\end{array}$$
and the observation that the vertical tangent bundle of the fibration 
$\pi$ is nothing else than $\cQ(1)$, we deduce that 
$$\begin{array}{rcl}
c_1(F(Y)) & = & 2h-3\s_1, \\
c_2(F(Y)) & = & -3h^2-9h\s_1+13\s_{11}-\s_2.
\end{array}$$

\begin{coro}\label{chern}
The Chern numbers and the arithmetical genus of the Fano surface $F(Y)$ are 
$$c_1(F(Y))^2=23355, \qquad c_2(F(Y))=11961, \qquad 
p_a(F(Y))=2942.$$
\end{coro}

\noindent {\it Remark}. 
 We can compare with the invariants of the Fano surface $F_2(Z)$
of planes in a cubic fivefold $Z\subset\PP^6$. As a subvariety of 
$G(3,7)$, $F_2(Z)$ is defined as the zero locus of the global section of 
the rank ten vector bundle $S^3T^*$, defined by an equation of $Z$. 
Therefore 
$$[F_2(Z)]=c_{10}(S^3T^*)=1134\s_{442}+1701\s_{433},$$
as we can easily deduce from the expression given above for the Chern
classes of a third symmetric power. On $F_2(Z)$ we thus have 
$\s_2=1134$ and $\s_1^2=2835$. Since $c(TF)=1-3\s_1+13\s_{11}-\s_2$, 
we conclude that 
$$c_1(F_2(Z))^2=25515, \qquad c_2(F_2(Z))=13041, \qquad 
p_a(F_2(Z))=3212.$$
In particular we have the relation
 $p_a(F(Y))=p_a(F_2(Z))+1-g(\Gamma(Y))$,
as expected.

\smallskip
We can also deduce the number $\nu$ of conics in $Y$ passing through a general
point. For this consider the universal conic $\cC\ra\cH$. This is a
subscheme of the universal supporting plane $\cP=\PP(\pi^*T)
\stackrel{\rho}{\lra}\cH$. Let $\gamma=\pi\circ\rho$. 
On $\cP$ we have  a tautological line bundle
$\cO_{\cP}(-1)\subset\gamma^*T$, and a tautological section $\tau$
of $\rho^*\cO_{\cH}(1)\otimes\cO_{\cP}(2)$ defined by the composition 
$$\rho^*\cO_{\cH}(-1)\hookrightarrow \gamma^*S^2T^*\lra \cO_{\cP}(2);$$
the zero locus of $\tau$ is precisely $\cC$. 

Let $V\subset\CC^6$ be a general three dimensional subspace. The zero
locus $\cP_V$ of the induced morphism
$$\cO_{\cP}(-1)\subset\gamma^*T\lra (\CC^6/V)\otimes \cO_{\cP}$$
is the set of points in $\cP$ mapped to $\PP V$ by the natural
morphism $\cP\lra\PP^5$. In particular the class of $\cP_V$ is $H^3$, 
if $H$ denotes the first Chern class of $\cO_{\cP}(1)$, and the class 
of its intersection $\cC_V$ with $\cC$ is $(2H+h)H^3$. 

Now the intersection of the cycle $\cC_V$ with $\rho^{-1}(F(Y))$ is 
the number of conics in $Y$ meeting $\PP V$. Since $Y$ has degree 
6 this number is equal to $6\nu$, and therefore
$$6\nu=(2H+h)H^3\rho^*[F(Y)]=\rho_*(2H^4+hH^3)[F(Y)].$$
But $\rho_*H^3=-\s_1$ and  $\rho_*H^4=\s_{11}$, so using Proposition 
\ref{fund} we obtain:

\begin{prop}
There are $\nu=318$ conics in $Y$ passing through a given general point.  
\end{prop}

Along the same line of ideas, one can consider in $F(Y)$ the curve
$\Delta$ of degenerate conics, and the Steiner map $s :\Delta\ra\PP^5$
mapping each such conic to its vertex. The same arguments as in 
\cite{CMW} yield
$$\deg s(\Delta) = (2h+2\s_1)(3h+2\s_1)=4860.$$

\subsection*{3.2 The Abel-Jacobi isomorphism for $F(Y)$}

Since Clemens and Griffiths proved that the Abel-Jacobi mapping 
$Alb\,F(X)\ra J(X)$ is an isomorphism for $X$ a general cubic threefold
and $F(X)$ is Fano variety of lines, similar statements have been obtain
for several other Fano manifolds. We fill in a gap in the literature 
by proving the following statement.  

\begin{theo}\label{AJ}
Let $Y$ be a general prime Fano threefold of genus 4.
Then the Abel-Jacobi mapping 
$$Alb\, F(Y)\lra J(Y)$$
is an isomorphism; and so by duality $J(Y)=J(Y)^*\simeq 
Alb\, F(Y)^*=Pic^0F(Y).$
\end{theo}

We adopt the strategy developed by Clemens and Griffiths, 
and used in \cite{Let} for the quartic
threefold, in \cite{CV} for the sextic double solid and 
in \cite{Col} for the cubic fivefold. 
The claim is that the theorem will follow from the existence of a 
Lefschetz pencil $(Y_t)_{t\in\PP^1}$ such that 
\begin{enumerate}
\item $Y_t=Q\cap C_t$ is contained in a fixed smooth quadric $Q$,
\item $Y_t$ and $F(Y_t)$ are smooth for general $t$,
\item if $Y_t$ is smooth but not $F(Y_t)$, then the later only has 
isolated singularities, 
\item if $Y_t$ is singular, then the singular locus of $F(Y_t)$ is 
a smooth curve, along which $F(Y_t)$ has two smooth branches intersecting 
transversely. 
\end{enumerate}

Consider a generic quadric hypersurface $Q\subset\PP^5$, a generic 
pencil of cubics $C_t\subset \PP^5$, $t\in\PP^1$, and then the 
pencil of complete intersections $Y_t=Q\cap C_t$. Condition (2)
clearly holds, and (3) follows from a simple dimension count. 
We focus on condition (4). The pencil $(Y_t)$ meets the discriminant 
hypersurface $\Delta$ parameterizing singular complete intersections of type 
$(2,3)$ at a finite number of points. These points must belong to the 
dense open subset of $\Delta$ parameterizing complete intersections
with a single node. So we are reduced to proving that for a general 
nodal complete intersection $Y$, the Fano surface of conics $F(Y)$ 
has the type of singularities allowed by condition (4). 

This involves the following steps, which we only sketch since the 
arguments are close from those of \cite{Let,CV,Col,PB}. 
\begin{enumerate}
\item {\it A conic $q\subset Y$ not passing through the vertex $v$ of $Y$, 
defines a  smooth point of $F(Y)$}. To check this we need to
characterize, for a conic $q\subset Y_{reg}$ ($Y$ having 
arbitrary singularities), the fact that it defines a smooth point 
in $F(Y)$. This goes as follows. Choose coordinates in $\PP^5$ such that 
$q$ be defined by the equations $x_3=x_4=x_5=0$ and ${\bar q}(x_0,x_1,x_2)=0$.  
If $Y=Q\cap C$, write the equations of $Q$ and $C$ as 
$$\begin{array}{rcl}
{\bar Q}(x) & = & \a {\bar q}(x_0,x_1,x_2)+x_3l_3+x_4l_4+x_5l_5, \\
{\bar C}(x) & = & \l {\bar q}(x_0,x_1,x_2)+x_3q_3+x_4q_4+x_5q_5. 
\end{array}$$
Suppose that $\a\ne 0$, which means that the supporting plane of $q$ is
not contained in $Q$. We say in that case that $q$ is a {\it non
      isotropic conic}. 
Then there is a unique cubic hypersurface $\a C-\l
      Q$ containing $Y$ and this supporting plane. Replacing $C$ by this
      cubic we may suppose that $\l=0$. Then we consider the exact
      sequence of normal bundles $0\ra N_{q/X}\ra N_{q/\PP^5}\ra
      N_{X/\PP^5|q}\ra 0$. We have $ N_{q/\PP^5}=\cO(1)^{\op 3}\oplus
      \cO(2)$ and $ N_{X/\PP^5}=\cO(2)\oplus \cO(3)$. 
Moreover, $F(Y)$ is smooth at $q$ if and only if $h^0(N_{q/X})=2$, 
which is equivalent to the surjectivity of the map 
  $H^0(N_{q/\PP^5})\ra H^0(N_{X/\PP^5|q})$. This map is easy to identify
      in terms of $\a,l_3,l_4,l_5,q_3,q_4,q_5$, and we come to the 
following conclusion: $q$ defines a smooth point of $F(Y)$ if and only
      if $q_3,q_4,q_5$ define linearly independent sections of
      $\cO(2)_{|q}$. 

Now suppose that $\a\ne 0$, that is $q$ is an {\it isotropic conic}.
We checked that the map  $H^0(N_{q/\PP^5})\ra H^0(N_{X/\PP^5|q})$
is always surjective under the hypothesis that $q\subset Y_{reg}$. 

\item {\it For $Y$ nodal but general, $F(Y)$ is smooth at any point 
defined by some conic $q\subset Y_{reg}$}. Indeed we can call a 
conic in some (arbitrary) $Y$ {\it special}, either if it is non 
isotropic but does not verify the condition above, or if it is 
isotropic and meets the singular locus of $Y$. Then the space of 
complete intersections $Y$ containing a special conic has at most 
two irreducible components, none of which coincides with the
      discriminant hypersurface.

\item {\it The set of conics in $Y$ passing through the vertex $v$, 
is the union of a complete curve $D$ parameterizing conics which are 
all smooth at $v$, and of $66$ conics $q(l,l')=l+l'$, where $l$ and $l'$
      are lines in $Y$ passing through $v$. Moreover 
the conics  $q(l,l')$ define
smooth points of $F(Y)$.} Indeed $Y$ contains $12$ lines passing through 
$v$; hence the $66$ reducible conics. The fact that for $Y$ general 
(among nodal complete intersections), a conic of type $q(l,l')$ 
defines a smooth point of $F(Y)$ is a boring computation. This
      computation also shows that the (linear) condition for a conic 
to pass through $v$ is transverse to the tangent space to $F(Y)$ at 
$q(l,l')$, which is thus isolated among the conics in $Y$ passing
      through $v$. 

\item {\it Let $\PP^+\ra\PP^5$ be the blow-up of $v$, and
      $Y^+\ra Y$ the strict transform of $Y$.
Let $F(Y^+)\stackrel{\pi}{\ra} F(Y)$ be the Fano surface of conics in $Y^+$.
Then $\pi^{-1}(D)=D_+\cup D_-$ is the union of two curves, and the
      restriction of $\pi$ to $D_{\pm}$ is a bijection with $D$. } 
More precisely, $F(Y^+)$ is the component of the Hilbert scheme of
      $Y^+$ containing the strict transforms of the conics in $Y$ not
      passing through $v$. If ${\tilde q}\in\pi^{-1}(q)$, where $q$ is a
      conic passing through $v$ but smooth at that point, then ${\tilde
      q}$ must be the union of the strict transform ${\bar q}$ of $q$, 
with a line $\ell$ in the exceptional divisor. Moreover $\ell$ must meet 
${\bar q}$ at the point defined by the tangent line to $q$ at $v$. 
But $S=E\cap Y^+$ is a smooth quadric surface, hence through that point 
pass exactly two lines $\ell_+$ and $\ell_-$, one from each of the two
      rulings of $S$. This is why $\pi^{-1}(q)=\{{\bar q}+\ell_+,
{\bar q}+\ell_-\}$, and $\pi^{-1}(D)=D_+\cup D_-$. 

\item {\it The Fano surface $F(Y^+)$ of conics in $Y^+$ is smooth along 
$D_{\pm}$, and the differential of the birational morphism 
$F(Y^+)\ra F(Y)$ maps the tangent
      planes to $F(Y^+)$ along $D_{\pm}$ to planes in the Zariski
      tangent space to $F(Y)$, meeting exactly along the tangent space
      to $D$.} The smoothness assertion follows from a computation 
with normal bundles, using the fact that $Y^+$ is a complete intersection in 
$\PP^+\subset\PP^5\times\PP^4$, of type $(1,1),(1,2)$. Moreover the
conic ${\bar q}+\ell_{\pm}$ is also a complete intersection, of type 
$(1,0),(0,1),(0,1),(1,1)$. If $p={\bar q}\cap\ell_{\pm}$, we  have
a commutative diagram
$$\begin{array}{ccccccccc}
 & 0 &  & 0 &  & 0 &  & 0 & \\
 & \downarrow &  & \downarrow &  & \downarrow &  & \downarrow & \\
 0\lra & H^0(N) & \lra & H^0(N(\ell_{\pm})) & \oplus & H^0(N({\bar q}))
 & \lra & N_p & \lra 0 \\
 & \downarrow &  & \downarrow &  & \downarrow &  & \downarrow & \\
 0\lra & H^0(E) & \lra & H^0(E(\ell_{\pm})) & \oplus & H^0(E({\bar q}))
 & \lra & E_p & \lra 0 \\
 & \downarrow &  & \downarrow &  & \downarrow &  & \downarrow & \\
 0\lra & H^0(F) & \lra & H^0(F(\ell_{\pm})) & \oplus & H^0(F({\bar q}))
 & \lra & F_p & \lra 0 \\
\end{array}$$ 
where $N$ denotes the normal bundle of ${\tilde q}$ in $Y^+$ (recall
      that $Y^+$ is smooth and ${\tilde q}$ is a locally complete
      intersection), $E$ denotes the normal bundle of ${\tilde q}$ in
      $\PP^+$, and $F$ the normal bundle of $Y^+$  in
      $\PP^+$, restricted to ${\tilde q}$. Moreover $N(\ell_{\pm})$ denotes
      the restriction of $N$ to $\ell_{\pm}$, and so on. 

Denote by $\phi(\ell_{\pm})$ the map $H^0(E(\ell_{\pm}))\ra 
H^0(F(\ell_{\pm}))$. The smoothness of the tangent cone to $C$ at the
      vertex $v$ is enough to imply that $\phi(\ell_{\pm})$ is
      surjective, hence $h^0(N(\ell_{\pm}))=2$. Moreover there is no 
non trivial section vanishing at $p$, and the evaluation map 
$H^0(N(\ell_{\pm}))\ra N_p$ is an isomorphism. Therefore 
$H^0(N)\simeq H^0(N({\bar q}))$. 

Now denote by  $\phi({\bar q})$ the map $H^0(E({\bar q}))\ra 
H^0(F({\bar q}))$. The fact that $q$ is non special implies that 
$\phi({\bar q})$ is surjective. So $h^0(N)=h^0(N({\bar q}))=2$
and $F(Y^+)$ is smooth at ${\tilde q}_{\pm}$. 

All this also holds for the finite number of isotropic conics 
passing through the vertex. Indeed one can check that if it did not
      hold for such an isotropic conic, then the conic would have to
      meet another singular point of $Y$, which is impossible.

Finally, we check that the 
images of the differentials of $\pi$ at   ${\bar q}+\ell_{\pm}$ 
meet along a line, which must be the Zariski tangent space to $D$. 
So $D$ is smooth, and isomorphic with $D_{\pm}$.   
\end{enumerate} 
This concludes the proof. \qed

\section{Application: Two integrable systems}

\subsection*{4.1 Cubic fivefolds containing a given fourfold}

By a result of Z. Ran, deformations of smooth Lagrangian subvarieties
are unobstructed \cite{Ran}. In particular the component $\cH$ of 
the Hilbert scheme of $F(X)$ containing  $F_2(Z)$, for $Z$ a 
general cubic fivefold containing $X$ as  hyperplane section, 
is generically smooth. Consider the universal family $\cF\ra\cH$. 
Donagi and Markman proved  that in such a
situation, the relative Picard bundle 
$Pic^0\cF\ra\cH$, is an algebraic completely integrable Hamiltonian
system (ACIHS), at least over the open subset $\cH_0\subset\cH$
parameterizing smooth deformations (\cite{DM1}, Theorem 8.1). 
This means that the total space of
the fibration admits a symplectic structure such that the fibers 
are Lagrangian subvarieties. 
Note that it follows  in particular that the dimension of $\cH$ 
equals $h^{1,0}(F_2(Z))=21$.  

The family $\cA$ of cubic fivefolds $Z$ containing $X=Z\cap H$ as a 
hyperplane section is a parameterized by the linear system 
$$|I_X(3)|=\langle c_X,x_0|\cO(2)|\rangle ,$$
where $c_X$ is an equation of $X$ in $H$ and $x_0$ is an equation of
$H$. This linear system has dimension $28$ and admits a natural action
of the seven-dimensional group 
$G$, consisting of automorphisms of $\PP^6$ whose restriction to 
$H$ is trivial. Moreover the subvariety 
$j_Z(F_2(Z))\subset F(X)$ does not change when $Z$ moves 
in a $G$-orbit. If we denote by $|I_X(3)|^{\circ}$ the open subset 
of $|I_X(3)|$ parameterizing smooth cubic fivefolds $Z$ such that 
$F_2(Z)$ is a smooth surface and $i_Z$ embeds  $F_2(Z)$ in $F(X)$, 
then the induced map 
$$|I_X(3)|^{\circ}\lra\cH$$
is constant on the $G$-orbits, and we may therefore consider
$\cH$ as a substitute for the moduli space of isomorphism classes 
of cubic fivefolds containing $X$. 

\medskip\noindent {\it Remark}.
Note that the tangent space to the corresponding deformation is 
$H^1(Z,TZ\otimes I_X)=H^1(Z,TZ(-1))$, which can be computed from 
the exact sequence
$$0\lra TZ(-1)\lra T\PP^6(-1)_{|Z}\lra \cO_Z(2)\lra 0.$$
The Euler exact sequence restricted to $Z$ gives
$H^0(Z,T\PP^6(-1)_{|Z})=H^0(\PP^6,\cO(1))^*:=V$ and 
$H^1(Z,T\PP^6(-1)_{|Z})=0$. Since $H^0(Z,\cO_Z(2))=H^0(\PP^6,\cO(2))
=S^2V^*$, we get that 
$$H^1(Z,TZ\otimes I_X)=Coker(V\stackrel{c_Z}{\lra} S^2V^*)$$
is the cokernel of the map $c_Z$ defined as the differential of an
equation of $Z$. This map is always injective for $Z$ smooth, hence 
$h^1(Z,TZ\otimes I_X)=21$. In particular we can easily obtain 
locally complete families for our deformation problem, parameterized by 
linear subsystems of $|I_X(3)|$ transverse to $\PP c_X(V)$. 

\medskip
It was proved by Collino that the Abel-Jacobi morphism 
$$Alb F_2(Z)\lra J(Z)$$
to the intermediate jacobian, is an isomorphism for a general 
cubic fivefold $Z$ \cite{Col}. 
In particular $Alb\,F_2(Z)$, like $J(Z)$, is self-dual, hence 
naturally isomorphic with the Picard variety $Pic^0F_2(Z)$. 
Denote by $\cH_s\subset\cH$ the open subset parameterizing the 
Lagrangian surfaces $i_Z(F_2(Z))\subset F(X)$, for $Z$ smooth 
with $F_2(Z)$ smooth. Denote by $\cH_a\subset\cH_s$
the open subset over which the Abel-Jacobi theorem does hold. 
We conclude:

\begin{theo}\label{acihs1}
The relative intermediate jacobian $J(\cH_a)$
over the family $\cH_a$ of smooth cubic fivefolds containing $X$, 
is an algebraic completely integrable Hamiltonian system. 
\end{theo}

In particular, the {\it cubic condition} of Donagi and Markman 
must hold (see \cite{DM1}, section 7.2). That is, the map 
$$T_{J(Z)}\cA=Sym^2H^{3}(Z,\Omega_Z^2)\lra T_{[Z]}\cH_a\simeq
H^1(Z,TZ\otimes I_X),$$
where $\cA$ denotes the moduli space of polarized abelian
varieties,
must be completely symmetric. To be precise, note that 
$H^{3}(Z,\Omega_Z^2)$ can easily be computed with the help of the normal
exact sequence for $Z\subset\PP^6$ and the Euler sequence restricted to
$Z$. We obtain 
$$H^{3}(Z,\Omega_Z^2)=Ker(S^2V\stackrel{c_Z^t}{\lra}V^*)
=H^1(Z,TZ\otimes I_X)^*.$$
Let $H:=H^{3}(Z,\Omega_Z^2)$. The cubic condition is that the map
$S^2H\ra H^*$ introduced above, which is just the differential of 
the induced map from $\cH_a$ to the moduli space $\cA$,
 comes from a cubic form $\theta_Z\in
S^3H^*$. This form can be described as the tautological cubic 
$$\theta_Z: S^3H\hookrightarrow S^3(S^2V)\lra S^2(S^3V)
\stackrel{c_Z^2}{\lra}\CC.$$

\subsection*{4.2 Fano threefolds of genus four  contained 
in  a given cubic fourfold}

We want to extend our ACIHS  to an open subset
$\cH_n\subset\cH$, containing $\cH_s$,  parameterizing  
Lagrangian surfaces $i_Z(F_2(Z))\subset F(X)$, where $Z$ is allowed 
to be a general cubic fivefold with a node. In this case we know that 
$F_2(Z)$ gets singular along a curve, and that its normalization 
is the Fano surface of conics $F(Y)$ of a general prime Fano threefold 
of genus 4.  
 
We begin with a few observations. Let $Z$ be a general nodal 
cubic fivefold containing $X$. 
\begin{enumerate}
\item The singular surface $F_2(Z)$ is a flat deformation of the smooth 
$F_2(Z')$, with $Z'$ smooth. This follows e.g. 
from Kollar's
criterion (6.1.3) in \cite{Kol}, which applies since $F_2(Z)$ is
a locally complete intersection, hence has property  $S_2$.
\item The Hilbert scheme $\cH$ remains smooth at $i_Z(F_2(Z))$. This is
because $F_2(Z)$ is a locally complete intersection and $i_Z$ is a
closed embedding. So  $i_Z(F_2(Z))$ is again a locally complete 
intersection in $F(X)$, again Lagrangian, so Ran's result on the 
unobstructedness of deformations remains valid (see \cite{Kaw}, 
and also \cite{FM}).
\end{enumerate}

\begin{prop}
The ACIHS $Pic^0F_2(\cZ)\lra\cH_s$ extends to $\cH_n$.
\end{prop} 

\proof Over $\cH_s$ we have a natural identification between
$H^{0,1}(F_2(Z))$ (the underlying vector space of $Pic^0F_2(Z)$)
and the dual to the tangent space of the base, that is 
$H^0(N_{i_ZF_2(Z)/F(X)})$, which is isomorphic to $H^{1,0}(F_2(Z))$
since $i_ZF_2(Z)$ is Lagrangian. We need to prove that this 
identification extends to $\cH_n$. Then the claim follows since the
cubic condition of Donagi and Markman is closed, and 
also the property of the symplectic 
form of being closed.

First we describe, for $Z$ a general nodal cubic, the relevant data 
in terms of the associated Fano threefold $Y$ of genus 4. Recall that 
the map 
$\nu: F(Y)\ra F_2(Z)$ restricts to a birational 
isomorphism between $F(Y)-B_{\pm}(Y)\ra F_2(Z)-\Gamma(Y)$, and maps 
the curves $B_{+}(Y)$ and $B_{-}(Y)$ isomorphically to $\Gamma(Y)$. 
Moreover the two branches of $F_2(Z)$ along $\Gamma(Y)$ intersect
 transversely at every point, and this remains true for $i_ZF_2(Z)$
since $i_Z : F_2(Z)\ra F(X)$ is a closed embedding. We denote 
for simplicity $\Gamma=i_Z(\Gamma(Y))$ and $F=i_ZF_2(Z)$. 

Since $F$ is a locally complete intersection in $F(X)$, its normal
sheaf $N_F=Hom(\cI_F/\cI_F^2,\cO_F)$ is a locally free $\cO_F$-module
of rank two. Since the two branches of $F$ along $\Gamma$ are
Lagrangian, their tangent spaces generate in $TF(X)$ the orthogonal 
$T\Gamma^{\perp}$ to 
$T\Gamma$ with respect to the symplectic form. The quotient   
of $T\Gamma^{\perp}$ by $T\Gamma$ is the direct sum of the normal 
bundles $N_+$ and $N_-$ of $B_+(Y)$ and $B_-(Y)$ in $F(Y)$. In
particular the symplectic form restricts to a non degenerate pairing 
$N_+\otimes N_-\lra\cO_{\Gamma}$. 

Consider the tangent map
$TF(Y)\ra\nu^*TF(X)\simeq\nu^*\Omega^1_{F(X)}$. Dualizing and 
pushing forward, we get a map $TF(X)\otimes\cO_F\ra\nu_*\nu^*TF(X)
\ra\nu_*\Omega^1_{F(X)}$, whose image we denote by $N_0\subset
\nu_*\Omega_{F(Y)}^1$. Clearly $N_0$ and $\nu_*\Omega_{F(Y)}^1$
are isomorphic outside $\Gamma$. At a point $p\in\Gamma$, an element 
the fiber of $\nu_*\Omega_{F(Y)}^1$ consists in  pairs of differential 
forms on the two branches of $F_2(Z)$ at that point. Such a pair is 
in the image of $TF(X)\otimes\cO_F\simeq \Omega^1_{F(X)}\otimes\cO_F$
if and only if the restrictions of the two forms to $T_p\Gamma$
coincide. So the quotient of $\nu_*\Omega_{F(Y)}^1$ by $N_0$ is
naturally identified with $\Omega_{\Gamma}^1$. 

On the other hand, there is a natural map from $N_0$ to $N$, which is 
again an isomorphism outside $\Gamma$. At a point $p\in\Gamma$, we can
choose local coordinates $x,y,z,t$ on $F(X)$ such that $F$ is given
by the equations $xy=z=0$, and $\Gamma$ by $x=y=z=0$. An element in the 
fiber of $N_0$ at $p$ can be seen as a linear form on $\langle
dx,dy,dz\rangle$, and it is mapped in $N$ to the corresponding linear
form on $\langle d(xy)=xdy+ydx,dz\rangle$. This maps $N_0$ injectively in 
$N$, with image the subsheaf of $N$ defined locally by the condition
that the image of $xy$ vanishes along $\Gamma$. Globally, this means
that $N_0$ is the kernel of the natural map $N\ra N_+\otimes N_-$. 
And we have noticed that this line bundle on $\Gamma$ is trivialized 
by the symplectic form. Note that in terms of deformations, 
$N_0$ is the subsheaf 
of $N$ that preserves the singularity $xy=0$ at first order. 

We get the following diagram:

$$\begin{array}{ccccccc}
 & 0 && & & & \\
 &  \downarrow &&  & & & \\
0  \lra & N_0 & \lra & \nu_*\Omega_{F(Y)}^1 & \lra & \Omega_{\Gamma}^1
& \lra  0 \\
 &  \downarrow && & & & \\
 &  N && & & & \\
 &  \downarrow && & & & \\
 &  \cO_{\Gamma} && & & & \\
 &  \downarrow && & & & \\
 &  0 && & & & 
\end{array}$$ 
  
Taking global sections, we get that $h^0(N)\le 1+h^0(N_0)\le 
1+h^0(\Omega_{F(Y)}^1)$. But $h^0(N)=21$ and $h^{1,0}(F(Y))=20$, 
so that these inequalities are equalities, and we obtain the exact 
sequence
$$0\lra H^0(N_0)\simeq H^0(\Omega_{F(Y)}^1)\lra H^0(N)\lra 
H^0( \cO_{\Gamma})\lra 0.$$
This must be interpreted as follows: $H^0(N)$ is the tangent space to
the Hilbert scheme of $F(X)$ at the point defined by $F_2(Z)$ (recall
that the Hilbert scheme is smooth at that point); $H^0(N_0)\simeq 
H^0(\Omega_{F(Y)}^1)$ is the tangent space to the hypersurface
parameterizing the surfaces $i_ZF_2(Z')$ for $Z'$ a nodal cubic, or
equivalently the images of the Fano surfaces $F(Y)$ for $Y$ a Fano
threefold of genus 4 contained in $X$. And $H^0(\Omega_{F(Y)}^1)$ 
is naturally isomorphic with the first order deformations of $Y$ (see ?). 

On the other hand, we have an exact sequence 
$0\ra \cO_F\ra\nu^*\cO_{F(Y)}\ra\cO_{\Gamma}\ra 0,$ (where the map
to $\cO_{\Gamma}$ is given by the difference of the restrictions to 
$B_+(Y)$ and $B_-(Y)$). Hence the exact sequence 
$$0\lra H^0(\cO_{\Gamma})\lra H^1(\cO_F)\lra 
H^1( \cO_{F(Y)})\lra 0.$$

Now we consider the relative Picard fibration $Pic^0(\cF)\ra\cH_n$. 
The symplectic form defined over $\cH_s$ has a natural extension to 
$\cH_n$ which can be defined as follows (see \cite{DM1}, section 8.4). 
If $L$ is an invertible sheaf over $F_2(Z)$, considered as a torsion
sheaf on $F(X)$, then the  tangent space 
to $Pic^0(\cF)$ at $L$ is naturally identified with 
$Ext^1_{\cO_{F(X)}}(L,L)$. And there is a natural skew-symmetric form on
this space, defined by the composition 
$$\wedge^2Ext^1_{\cO_{F(X)}}(L,L)
\lra Ext^2_{\cO_{F(X)}}(L,L)\stackrel{c_1(\cO(1))}{\lra}
Ext^4_{\cO_{F(X)}}(L,L)\simeq Hom(L,L)^*=\CC,$$
where we have used Serre duality and the fact that the canonical sheaf 
of $F(X)$ is trivial. Now the local to global spectral sequence for 
$Ext$'s easily yields an exact sequence
$$0\lra H^1(\cO_F)\lra Ext^1_{\cO_{F(X)}}(L,L)\lra H^0(N)\lra 0,$$
which is nothing else than the tangent sequence of the Picard
fibration. Since the fibration is Lagrangian over $\cH_s$, by continuity
$H^1(\cO_F)$ is isotropic. So the skew-symmetric form on 
$Ext^1_{\cO_{F(X)}}(L,L)$ is non degenerate if and only if the induced
pairing 
$$H^1(\cO_F)\otimes H^0(N)\lra H^1(N)\stackrel{c_1(\cO(1))}{\lra}
H^2(N\otimes\Omega^1_{F(X)|F})\lra H^2(\det N)=H^2(\omega_{F})=\CC$$
is non degenerate. To state it in a more convenient form, we need to
prove that the map 
$$\CC\stackrel{c_1(\cO(1))}{\lra}H^1(\Omega^1_{F(X)|F})
\lra H^1(N)\lra Hom(H^1(\cO_F),H^0(N)^*)$$ 
maps $c_1(\cO(1))$ to an isomorphism. Over $\cH_s$, $N$ is identified
with $\Omega^1_F$ and this follows from the
Hard Lefschetz theorem. Over $\partial\cH_n:=\cH_n-\cH_s$ we use our exact 
sequences
$$\begin{array}{ccccccc}
0\lra & H^0(\cO_{\Gamma}) & \lra & H^1(\cO_F) & \lra & 
H^1(\cO_{F(Y)}) & \lra 0 \\
 & \downarrow & & \downarrow && \downarrow & \\
0\lra & H^1( \Omega^1_{\Gamma}) & \lra & H^0(N)^* & 
\lra & H^0(\Omega^1_{F(Y)})^* & \lra 0 
\end{array}$$
The restricted map $H^0(\cO_{\Gamma})\ra H^0(N)^*\ra
  H^0(\Omega^1_{F(Y)})^*$ is
 zero. This can be seen as follows: there is a commutative diagram
$$\begin{array}{ccccc}
H^0(N)\otimes H^1(\cO_F) & \lra & H^1(N) & \stackrel{c_1(\cO(1))}{\lra}
& H^2(\omega_F) \\
  \uparrow & &\uparrow & &\uparrow  \\
H^0(N_0)\otimes H^0(\cO_{\Gamma}) & \lra & H^1(N_0\otimes\cO_{\Gamma}) 
& \stackrel{c_1(\cO(1))}{\lra} & H^1(\omega_F\otimes\cO_{\Gamma}) 
\end{array}$$
where the right arrow of the bottom line involves the sheaf 
morphisms $$N_0\otimes\Omega_{F(X)}^1\ra N_0\wedge N\ra\det N=\omega_F.$$ 
But the image of $\Omega_{F(X)}^1$ in $N$ is precisely $N_0$, and 
$\wedge^2N_0\otimes\cO_{\Gamma}=0$. 

Now,  the induced map $$H^1(\cO_{F(Y)})\lra
H^0(\Omega^1_{F(Y)})^*\simeq H^2(\Omega^1_{F(Y)})$$ is defined 
by the cup product with $c_1(\nu^*\cO(1))$. Since $\nu$ is finite 
$\nu^*\cO(1)$ is still ample, so Hard Lefschetz applies and this 
map is an isomorphism. On the other hand the other induced map 
$H^0(\cO_{\Gamma})\lra H^1( \Omega^1_{\Gamma})$ is given by 
$c_1(\cO(1)_{|\Gamma})$, and again it is an isomorphism. This concludes
  the proof.  \qed

\begin{prop} 
Let $Z$ be a general nodal cubic and $Y$ the associated 
prime Fano threefold of genus 4.  
There is an exact sequence 
$$0\lra \CC^*\lra Pic^0F_2(Z)\stackrel{\nu^*}{\lra} 
Pic^0F(Y)\simeq J(Y)\lra 0.$$
\end{prop}

\proof The fact that $Pic^0F(Y)\simeq J(Y)$ is Theorem \ref{AJ}. 
Now consider the exact sequence
$0\ra\cO_{F_2(Z)}^*\ra\nu_*\cO_{F(Y)}^*\ra\cO^*_{\Gamma(Y)}\ra 0$, 
and the associated long exact sequence
$$0\lra H^0(\cO^*_{\Gamma(Y)})=\CC^*\lra Pic\; F_2(Z)\stackrel{\nu^*}{\lra}
 Pic\; F(Y).$$
The pull-back of a Weil divisor from $F_2(Z)$ to $F(Y)$ has 
the same intersection number with the curves $B_+(Y)$ and 
$B_-(Y)$, and this property characterizes the image of $\nu^*$. 
In particular it contains $Pic^0F(Y)$. Finally, 
$H^0(\cO^*_{\Gamma(Y)})$ being connected obviously maps to $Pic^0F_2(Z)$. 
Hence the claim. \qed

\medskip Recall that the Fano threefold $Y$ has been defined as the
intersection in $H=\PP^5$ of the fixed cubic fourfold $X$, with the 
base $Q$ of the tangent cone to $Z$ at its unique node $p_Z$.
Since $G$ fixes $H$, the map $Z\mapsto Q$ is constant along the 
$G$-orbits. 

\begin{prop}
There is a birational isomorphism 
$$|\cO_{\PP^5}(2)|\lra \partial\cH_n=\cH_n-\cH_s.$$
\end{prop} 

\proof Let $\Delta_X\subset |I_X(3)|$ denote the discriminant 
hypersurface, parameterizing the singular cubics containing $X$. 
For $X$ general $\Delta_X$ is irreducible. Let $\Delta_X^n$ denote the
open subset parameterizing cubics with a single node. The map 
$\Delta_X^n\ra |\cO_{\PP^5}(2)|$ sending the nodal cubic $Z$
to the trace on $H$ of its tangent cone at the node, is constant on the 
$G$-orbits. More precisely, it is easy to check that its fibers are 
precisely the $G$-orbits. This implies that there exists a $G$-stable 
open subset 
$\Delta_X^0\subset\Delta_X^n$ such that the restriction map 
$\Delta_X^0\ra |\cO_{\PP^5}(2)|$ is a good quotient of the $G$-action. 
In particular, there is an induced rational map $|\cO_{\PP^5}(2)|
\ra\partial\cH_n$, which is clearly dominant. 

So we just need to check that this map is bijective, that is, 
a general Fano threefold $Y=X\cap Q$ is uniquely determined by 
the singular Lagrangian surface $i_Z(F_2(Z))\subset F(X)$ associated 
to the nodal cubic $Z$ whose equation is $P+x_0Q=0$ for some equation 
$x_0$ of $H\subset\PP^6$. By construction, the singular locus 
of $i_Z(F_2(Z))$ is nothing else than the curve $\Gamma(Y)\subset F(X)
\subset G(2,6)$ of lines in $Y$. So we just need to check that
$\Gamma(Y)$ defines uniquely the quadric $Q$ such that $Y=X\cap Q$. 
That is, we must prove that 
$$H^0(G(2,6),I_{\Gamma(Y)}\otimes S^2T^*)=\CC.$$
For $Y$ generic, $I_{\Gamma(Y)}$ can be resolved by the Koszul complex 
of the section of the vector bundle $E=S^2T^*\oplus S^3T^*$ that defines 
$\Gamma(Y)$. The claim then follows from the identities
$H^0(G(2,6),E^*\otimes S^2T^*)=\CC$ and $H^k(G(2,6),\wedge^{k+1}E^*
\otimes S^2T^*)=0$ for $k>0$. Both facts are easy consequences of 
Bott's theorem on the Grassmannian. \qed 

\medskip
Now we are exactly in the situation considered in \cite{DP}: an 
ACIHS $\cF\ra B$ is defined over some smooth base $\cB$, and there is a 
hypersurface $\Delta\subset\cB$ over which the fibers degenerate
to extensions 
$$0\lra \CC^*\lra \cF\lra J(\cF)\lra 0,$$ 
where $J(\cF)$ is a family of abelian varieties over $\Delta$. 
If $h$ is a local equation of $\Delta$, suppose that the 
corresponding Hamiltonian vector field, which is a vertical vector field
for the fibration,  is tangent to the $\CC^*$ direction.   
Then symplectic reduction applies, and there is an induced ACIHS 
$J(\cF)\ra\Delta$. 

We apply this result to our setting: $\cF\ra B$ is our relative Picard
fibration $Pic^0(F_2(\cZ))\ra\cH_n$. Over the hypersurface $\Delta=
\partial\cH_n$ defined by nodal cubics, the fibers $Pic^0(F_2(\cZ))$
are $\CC^*$-extensions of the intermediate jacobians $J(Y)$ of the
associated Fano threefolds $Y$ of genus 4. Moreover we have seen that 
the $\CC^*$ factor is orthogonal, with respect to the extended
symplectic form, to the tangent hyperplane to $\Delta$. So symplectic 
reduction applies and we conclude: 

\begin{theo}\label{acihs2}
Let $X\subset\PP^5$ be a general cubic fourfold.  
Consider inside the 20-dimensional linear system $|\cO_{\PP^5}(2)|$, 
the open subset $|\cO_{\PP^5}(2)|_X$ of smooth quadrics $Q$ transverse to $X$. 
Denote by $$J(\cY_X)\lra |\cO_{\PP^5}(2)|_X$$ the family of intermediate 
jacobians of the genus 4 prime Fano threefolds $Y=X\cap Q$, where 
$Q\in |\cO_{\PP^5}(2)|_X$. 
Then $J(\cY_X)$ is an algebraic completely
 integrable Hamiltonian system. 
\end{theo}

Here again the cubic form can easily be identified. Let 
$W=H^0(\cO_{\PP^5}(1))^*$. From the normal 
sequence of $Y$ we can easily deduce an exact sequence 
$$H^0(Y, T\PP^6(-1)_{|Y})=W\stackrel{c_Q,c_X}{\lra}
W^*\oplus S^2W^*/Q\lra 
H^1(Y,TY(-1))\simeq H^1(Y,\Omega_Y^2)\lra 0.$$
Since $Q$ smooth, $c_Q$ is an isomorphism and we get an 
isomorphism $S^2W^*/Q\simeq
H^1(Y,\Omega_Y^2)$. The dual space $H^1(Y,\Omega_Y^2)$
can therefore be indentified with the with the space $A_Q\subset S^2W$
corresponding to quadrics in the dual projective space which are apolar
to $Q$. The cubic form $\theta_Q$ on $A_Q$ is then simply given by 
the composition
$$\theta_Q: S^3A_Q\hookrightarrow S^3(S^2W)\lra S^2(S^3W)
\stackrel{c_X^2}{\lra}\CC.$$


\vspace{1cm}

\noindent 
{\bf Atanas Iliev}\\ 
Institute of Mathematics,  
Bulgarian Academy of Sciences,  
Acad. G. Bonchev Str., bl. 8\\     
1113 Sofia, Bulgaria\\  
{\bf e-mail:} ailiev@math.bas.bg 

\bigskip

\noindent
{\bf Laurent Manivel}\\  
Institut Fourier, 
Laboratoire de Math\'ematiques, 
UMR 5582 (UJF-CNRS), BP 74\\    
38402 St Martin d'H\`eres Cedex, France\\ 
{\bf e-mail:}  Laurent.Manivel@ujf-grenoble.fr

\end{document}